\documentclass[12pt]{amsart}
\usepackage{}
\usepackage{stmaryrd}
\usepackage{amsfonts}
\usepackage{bbm}
\usepackage{amsmath}
\usepackage{amssymb}
\usepackage{a4, latexsym, amsfonts}
\usepackage[colorlinks=true]{hyperref}

\hypersetup{colorlinks=true,CJKbookmarks=true}

 \newtheorem{theorem}{Theorem}[section]
 \newtheorem{lemma}[theorem]{Lemma}

 \newtheorem{conjecture}[theorem]{Conjecture}
 
 \newtheorem{remark}[theorem]{Remark}
 \numberwithin{equation}{section}

\begin{document}

\title[]{Automorphism of reduced finitary power monoids}
\author[]{Yongjie Peng$^{*}$}
\author[]{Guixin Deng}

\date{}
\maketitle
 \edef \tmp {\the \catcode`@}
   \catcode`@=11
   \def \@thefnmark {}


\begin{abstract}Let \(H\) be a multiplicatively written monoid with identity element \(1_H\). The family \(\mathcal{P}_{\mathrm{fin},1}(H)\) of all finite subsets of \(H\) containing \(1_H\) is a monoid with identity $\{1_H\}$ under the setwise multiplication induced by \(H\). We prove that \(\operatorname{Aut}(\mathcal{P}_{\mathrm{fin},1}(H)) \cong \operatorname{Aut}(H)\) whenever  \(H\) satisfies either of the following two conditions:
(1) \(H\) is commutative, cancellative, and contains an invertible element of infinite order;
(2) \(H\) is an abelian  group which is not isomorphic to the Klein four-group.


\end{abstract}

\section{Introduction}
\label{sec:intro}

Let \(H\) be a semigroup. Endowed with the setwise multiplication induced by \(H\), the family of all nonempty subsets of \(H\) forms a semigroup, hereinafter denoted by \(\mathcal{P}(H)\).
The family \(\mathcal P_{\mathrm{fin}}(H)\) of all nonempty finite subsets of \(H\) is a subsemigroup of \(\mathcal{P}(H)\). More explicitly, the operation is defined by
\[
    XY := \{xy : x \in X,\; y \in Y\},
\qquad\text{for all } X, Y \in \mathcal P(H).
\]
We call \(\mathcal P(H)\) and \(\mathcal P_{\mathrm{fin}}(H)\) the large power semigroup and the finitary power semigroup of \(H\), respectively.

In particular, if \(H\) is a monoid with identity \(1_H\) then $\mathcal{P}(H)$ and $\mathcal P_{\mathrm{fin}}(H)$ are monoids too; their identity is the singleton $\{1_H\}$. Moreover,
\[
\mathcal P_1(H) := \{X \subseteq H : 1_H \in X\}
\quad\text{and}\quad
\mathcal P_{\mathrm{fin},1}(H) := \mathcal P_1(H) \cap \mathcal P_{\mathrm{fin}}(H)
\]
 are submonoids of \(\mathcal{P}(H)\) and \(\mathcal{P}_{\mathrm{fin}}(H)\), respectively. We call $\mathcal{P}_1(H)$ the reduced large power monoid and $\mathcal{P}_{\mathrm{fin},1}(H)$ the reduced finitary power monoid of \(H\).

\begin{remark}
{\textup{When}} $H$ {\textup{is an additively written semigroup, the same convention is adopted for}} $\mathcal{P}(H)$ {\textup{and its subsemigroups. The operation on}} $\mathcal{P}(H)$ {\textup{is then called setwise addition, and it is explicitly defined by}}
\[
X+Y := \{x+y: x \in X, y \in Y\},
\qquad \textup{for all}~X, Y \in \mathcal{P}(H).
\]
{\textup{Accordingly, if}} $H$ {\textup{is a monoid, we use}} $\mathcal{P}_0(H)$ {\textup{in place of}} $\mathcal{P}_1(H)$ {\textup{and}} $\mathcal{P}_{\mathrm{fin},0}(H)$ {\textup{in place of}} $\mathcal{P}_{\mathrm{fin},1}(H)$.
\end{remark}

Power semigroups first appeared explicitly in a 1950 paper of Ballieu \cite{bal-1950}, and their systematic study was started by Tamura and Shafer \cite{Shaf-1967,Tam-Sha1967} in the late 1960s. In recent years, they have been studied from several perspectives, including their algebraic structure, their factorization properties, and their connections with arithmetic combinatorics;  see \cite{An-Tr18,Cos-Tri2025,Fa-Tr18}. For more details about the arithmetic of power monoids, we refer the reader to the survey \cite{Tri-2026(a)}.

Let \(\mathbb{N}\) and $\mathbb{Z}$ be the additive monoid of  nonnegative integers and  the additive group of integers, respectively. Denote by $\mathrm{Aut}(H)$ the automorphism group of a semigroup $H$. The study of the automorphisms of power semigroups can be traced back to Byrd et al. In particular, it is proved in \cite[Theorem 4]{Byrd78} that  $\mathrm{Aut}(\mathcal{P}_{\mathrm{fin}}(H))$ is isomorphic to $\mathrm{Aut}(H)$ for a locally finite group $H$ with $|H|>5$,  and in \cite[Theorem 9]{Byrd84}  that  $\mathrm{Aut}(\mathcal{P}_{\mathrm{fin}}(\mathbb{Z}))$ is a splitting extension of $\mathbb{Z}$ by the Klein four-group.  The latter result does not completely determine the structure of $\mathrm{Aut}(\mathcal{P}_{\mathrm{fin}}(\mathbb{Z}))$. This was recently achieved  by Tringali and Wen \cite[Theorem 3.5]{Tri-2025}, who proved that   $\mathrm{Aut}(\mathcal{P}_{\mathrm{fin}}(\mathbb{Z}))$ is isomorphic to the direct product of a cyclic group of order two and the infinite dihedral group.  Wong et al. \cite{Wong25-2} showed that for a numerical semigroup $H$ (that is, a cofinite subsemigroup of the nonnegative integers under addition), $\mathrm{Aut}(\mathcal{P}_{\mathrm{fin}}(H))$ is a cyclic group of order $2$ if $H$ is a discrete interval; otherwise $\mathrm{Aut}(\mathcal{P}_{\mathrm{fin}}(H))$ is trivial. Tringali and Wen \cite{Tri26(b)} proved that $\mathrm{Aut}(\mathcal{P}(H))$ is trivial for every numerical semigroup $H$, and that $\mathrm{Aut}(\mathcal{P}_0(H))$ is trivial whenever $0\in H$.

In this paper, we focus on the relation between the automorphism group of a monoid $H$ and  the automorphism group of the reduced finitary power monoid $\mathcal{P}_{\mathrm{fin},1}(H)$. Each $f\in \mathrm{Aut}(H)$ induces an automorphism $f^*$ of $\mathcal{P}_{\mathrm{fin},1}(H)$ defined as follows:
\[
    f^*(X)=f[X]:=\{f(x):x\in X\},
\qquad\text{for all }X\in\mathcal{P}_{\mathrm{fin},1}(H).
\]
We call  $f^*$  the augmentation of $f$ and we say that an automorphism of $\mathcal{P}_{\mathrm{fin},1}(H)$ is inner if it is the augmentation of an automorphism of $H$.

It is easy to verify that this gives a group monomorphism
\[
    \Phi:\operatorname{Aut}(H)\longrightarrow
    \operatorname{Aut}(\mathcal P_{\mathrm{fin},1}(H)),
    \qquad f\longmapsto f^*.
\]

It is natural to ask   whether $\Phi$ is an isomorphism. In other words, is each automorphism of $\mathcal{P}_{\mathrm{fin},1}(H)$ inner?   Tringali and Yan \cite{Tri-Yan2025(b)} showed that  the only non-trivial automorphism of $\mathcal{P}_{\mathrm{fin},0}(\mathbb{N})$ is given by $X\longmapsto \max X-X$. This proves that, in general, $\Phi$ is not surjective since $\mathrm{Aut}(\mathbb{N})$ is trivial.  In fact, the same result was also obtained in \cite[Theorem 1]{Byrd84} but this paper went unnoticed until April 2026. More recently, Bhowmik and Tringali \cite[Theorem 1.2]{Bho-Tri26} showed that $\mathrm{Aut}(\mathcal{P}_{\mathrm{fin},0}(H))$ is trivial for every  numerical monoid $H\not=\mathbb{N}$, which confirmed a conjecture proposed in \cite{Tri-Yan2025(b)}.

Wong et al.~\cite[Theorem 1.3]{Wong25} proved that all automorphisms of $\mathcal{P}_{\mathrm{fin},0}(\mathbb{Z})$ are inner and hence gave an affirmative answer to a conjecture proposed by Tringali and Wen \cite{Tri-2025}. In a  remarkable feat, Rago \cite[Theorem 13]{Rago25(b)} established that if $H$ is a finite abelian group, then all automorphisms of $\mathcal{P}_{\mathrm{fin},1}(H)$ are inner, except when $H$ is the Klein four-group, in which case $\mathrm{Aut}(\mathcal{P}_{\mathrm{fin},1}(H))$ is isomorphic to the direct product of two copies of the symmetric group of degree three.

In this paper, we show that all automorphisms of $\mathcal{P}_{\mathrm{fin},1}(H)$ are inner if $H$ is a commutative and cancellative monoid containing an invertible element of infinite order or $H$ is an infinite abelian group.

\section{Notation and preliminary results}
\label{sec:2}
We use the conventions of \cite{Ho95} for semigroups and monoids. Let \(H\) be a monoid written   multiplicatively.  Then an element $x \in H$ is
right (resp. left) cancellative if  $yx=zx\Longrightarrow y=z$ (resp. $xy=xz\Longrightarrow y=z$) for any $y,z\in H$.
We call $H$ cancellative if each of its elements is both right and left cancellative. Let $X$ be a nonempty subset of $H$ and let $a\in H$. We write
$$aX:=\{ax:x\in X\}.$$
We denote by \(C_n\) the cyclic group of order \(n\). The cardinality of a set \(X\) is denoted by \(|X|\).

We start with a result of Tringali and Yan \cite{Tri26}.

\begin{lemma}\label{lem:pullback}\cite[Corollary 3.3]{Tri26}
Let \(H\) and \(K\) be monoids, and let \(F:\mathcal P_{\mathrm{fin},1}(H)\to \mathcal P_{\mathrm{fin},1}(K)\) be an isomorphism. Then there is a unique bijection \(f:H\to K\), called the pullback of \(F\), such that
\[
F(\{1_H,x\})=\{1_K,f(x)\},\qquad \text{for all } x \in H.
\]

In particular, $f(1_H)=1_K.$
\end{lemma}

The pullback $f$ of the isomorphism $F$ has many interesting properties and is remarkably  helpful  in the study of  reduced finitary power monoids.  In particular, we will need the following result by Tringali and Yan \cite{Tri26}.

\begin{theorem}\label{bothgroup}\cite[Theorem 5.1]{Tri26} Let $H$ and $K$ be cancellative monoids, and suppose at least one of them is torsion. Let $F$ be an isomorphism from $\mathcal P_{\mathrm{fin},1}(H)$ to $\mathcal P_{\mathrm{fin},1}(K)$. Then the pullback of $F$ is an isomorphism from $H$ to $K$, and both $H$ and $K$ are groups.
\end{theorem}

Rago \cite{Rago25(b),Rago26(c)} continued the study of the properties of the pullback  and proved the following results.

\begin{theorem}\label{lemma:2.4}\cite[Theorem 13]{Rago25(b)}
Let \(H\) be a finite abelian group with \(H\not\cong C_2\times C_2\). Then all automorphisms of $\mathcal{P}_{\mathrm{fin},1}(H)$ are inner.
\end{theorem}

\begin{theorem}\label{lemma:2.2}\cite[Theorem 14]{Rago26(c)}
Let \(H\) and \(K\) be commutative and cancellative monoids, each containing an element of infinite order. Suppose that at least one of them has a nontrivial unit group. Then the pullback of every isomorphism
\(F:\mathcal P_{\mathrm{fin},1}(H)\to\mathcal P_{\mathrm{fin},1}(K)\) is an isomorphism.
\end{theorem}

Let \(H\) be a finite abelian group with \(H\not\cong C_2\times C_2\). Indeed, Rago \cite{Rago25(b)} showed that every \(F\in\operatorname{Aut}(\mathcal P_{\mathrm{fin},1}(H))\)
whose pullback is the identity map is itself the identity map, and Theorem \ref{lemma:2.4} comes from this fact and Theorem \ref{bothgroup}. We summarize it as follows.

\begin{lemma}
\label{thm:2.5}Let $H$ be a monoid. Assume that the pullback of every automorphism of
\(\mathcal P_{\mathrm{fin},1}(H)\) is an automorphism of \(H\). Then the
following are equivalent:

$(1)$~ \(\Phi:\operatorname{Aut}(H)\longrightarrow
\operatorname{Aut}(\mathcal P_{\mathrm{fin},1}(H))\),
\(f\longmapsto f^*\), is an isomorphism;

$(2)$~ every \(F\in\operatorname{Aut}(\mathcal P_{\mathrm{fin},1}(H))\)
whose pullback is the identity map is itself the identity map;

$(3)$~if the pullbacks of $F_1,F_2\in \operatorname{Aut}(\mathcal P_{\mathrm{fin},1}(H))$ coincide, then $F_1=F_2$;

$(4)$~if $f$ is the pullback of an automorphism $F$ of $\mathcal P_{\mathrm{fin},1}(H)$, then $F(X)=\{f(x):x\in X\}$ for any $X\in \mathcal P_{\mathrm{fin},1}(H)$.

\begin{proof}
$(1)\Longrightarrow(2)$ is obvious.

$(2)\Longrightarrow (3)$. Let $f$ be pullback of $F_1$. By hypothesis, $f$ is also the pullback of $F_2$. On the other hand, it is clear that $f^{-1}$ is the pullback of $F_1^{-1}$, and hence $f^{-1} \circ f$  is the pullback of $F_1^{-1}\circ F_2$. So  $F_1^{-1}\circ F_2$ is the identity map and $F_1=F_2.$

$(3)\Longrightarrow (4)$. By assumption, $f$ is also an automorphism of $H$. So $F$ and $f^*$ share the same pullback $f$. Hence, $F(X)=f^*(X)=\{f(x):x\in X\}$ for any  $X\in \mathcal P_{\mathrm{fin},1}(H)$.

$(4)\Longrightarrow (1)$ is clear.
\end{proof}
\end{lemma}

We also need the following results about the reduced finitary power monoid of a cancellative monoid containing an element of infinite order.  The results were first proved by Rago in the commutative setting, see the proof of  \cite[Lemma 4]{Rago26(c)} and then generalized to the non-commutative setting by Bhowmik and Tringali \cite[Lemma 2.2]{Bho-Tri26}.

\begin{lemma}\label{lem:3.2}
Let \(H\) be a  cancellative monoid. Suppose  \(a\in H\) is invertible of infinite order and   $X$ is a finite nonempty subset of $H$.  Then there exists a  positive integer \(N\) such that
\[
a^{in}X\cap a^{jn}X=\emptyset
\]
for all \(n>N\) and \(i,j\in \mathbb{Z}\) with $i\neq j.$
\begin{proof} Lemma 2.2 in [4] gives that there exists a positive integer $N$ such that $a^{in} X \cap a^{jn} X =\emptyset$ for all integers $n>N$ and $i, j \in \mathbb {N}$ with $i \ne j$. Therefore, taking into account that $a$ is a unit of $H$, we get that $a^{(i+k)n} X$ and $a^{(j+k)n} X$ are also disjoint for all $i,j \in \mathbb N$ with $i\neq j$, $n>N$ and $k \in \mathbb Z$, which is the desired conclusion.
\end{proof}
\end{lemma}

\begin{lemma}\label{lemma:2.3}\cite[Proposition 2.4]{Bho-Tri26}
Let $H$ and $K$ be  cancellative monoids and assume that  $H$ contains  an element of infinite order. Then each isomorphism $F:\mathcal P_{\mathrm{fin},1}(H)\to\mathcal P_{\mathrm{fin},1}(K)$  is cardinality-preserving, that is, $|F(X)|=|X|$ for any
\(X\in\mathcal P_{\mathrm{fin},1}(H)\).
\end{lemma}

\section{Main Results}
\label{sec:3}

Let $H$ be a semigroup and let $h\in H$. The ``correlation function''
\[
\nu_h:\mathcal{P}_{\mathrm{fin}}(H)\longrightarrow\mathbb{N},\qquad A\longmapsto|\{g\in A: hg\in A\}|
\]
was first introduced in \cite[Definition 2.5]{Bho-Tri26}.
This function is closely related to the discrete covariogram in convex geometry \cite[Section 6.1]{Bianchi} and to difference representation functions in additive combinatorics \cite[Section 6.1]{Nathanson}. We need the following properties of $\nu_h$.

\begin{lemma}
\label{lemma:3.1}
Let \(H\) be a commutative and cancellative monoid. Then for any \(y,z\in H\) and
\(A\in\mathcal P_{\mathrm{fin}}(H)\),\\
$(1)$~ \(\nu_z(A)=|A\cap zA|\);\\
$(2)$~ \(\nu_z(A)=\nu_z(yA)\);\\
$(3)$~ if \(A=\bigcup_{i=1}^nA_i\) is a disjoint union, and
\( zA_i\cap A_j=\emptyset\) for all $i,j\in\{1,2,\ldots,n\}$ with $i\neq j$, then
$$\nu_z(A)=\sum_{i=1}^n\nu_z(A_i).$$
\begin{proof}
$(1)$~is a special case of  \cite[Lemma 2.6]{Bho-Tri26}.\\
$(2)$~Since $H$ is commutative and cancellative, the map
$$A\cap zA\longrightarrow yA\cap yzA:\qquad x\longmapsto xy,$$
is a bijection, so  $(2)$ follows from $(1)$.\\
(3) By the assumptions, we see that
$$A\cap zA=(\cup_{i=1}^nA_i)\cap(\cup_{j=1}^nzA_j)=\cup_{1\leq i,j\leq n}(A_i\cap zA_j)=\cup_{1\leq i\leq n}(A_i\cap zA_i).$$
Hence,
$$\nu_z(A)=|A\cap zA|=|\cup_{1\leq i\leq n}(A_i\cap zA_i)|=\sum_{i=1}^n|A_i\cap zA_i|=\sum_{i=1}^n\nu_z(A_i).$$
\end{proof}
\end{lemma}

\begin{lemma}\label{lem:3.3}\cite[Corollary 2.8]{Bho-Tri26}
Let $H$ be a commutative and cancellative monoid containing an element of infinite order, and let $F$ be an automorphism of $\mathcal{P}_{\mathrm{fin},1}(H)$ with identity pullback. Then $$\nu_z(A)=\nu_z(F(A))$$
for any $z\in H$ and $A\in \mathcal{P}_{\mathrm{fin},1}(H).$
\end{lemma}

Now we can state and prove our main results.

\begin{theorem}
\label{thm:infinite-unit}
Let \(H\) be a commutative and cancellative monoid, and suppose that $H$ contains  an invertible element of infinite order. Then all automorphisms of $\mathcal{P}_{\mathrm{fin},1}(H)$ are inner.

\begin{proof}Let \(P_h=\{1_{H},h\}\) for any \(h\in H\) and let $g\in H$ be an invertible element of infinite order.

By Theorem \ref{lemma:2.2}, the pullback of every isomorphism of $\mathcal{P}_{\mathrm{fin,1}}(H)$ is  an automorphism of $H$.
Thus, by Lemma~\ref{thm:2.5}, it suffices to show that if an automorphism $F$ of $\mathcal{P}_{\mathrm{fin},1}(H)$ has trivial pullback, then $F$ is the identity map. Accordingly, we assume below  that
\begin{equation}\label{3.1}F(P_h)=P_h\qquad\text{for every }h\in H.\end{equation}
Fix \(S\in\mathcal P_{\mathrm{fin},1}(H)\) and let \(T=F(S)\). We proceed to show that $S\subseteq T$. Take \(z\in S\).  Suppose $z \neq 1_H$, otherwise the conclusion is obvious.  The proof is divided into two cases.

Case 1: \(z\) is invertible. Since $U:=T\cup z^{-1}T$ is also finite, by  Lemma \ref{lem:3.2} there exists a positive integer $N$ such that   $g^{iN}U\cap g^{jN}U=\emptyset$ for any pairwise distinct $i,j\in\mathbb{Z}$.  Put $w=g^N$. Then we have
\begin{equation}\label{3.2}w^iR_1\cap w^jR_2=\emptyset,\qquad \text{for}~R_1,R_2\in\{T,z^{-1}T\}~\text{and}~i,j\in\mathbb{Z}~\text{with}~i\neq j.\end{equation}

Put $A=w^{-1}S\cup  wz^{-1}S\cup \{1_H\}$.  Since $1_H,z\in S$, we have $P_w\subseteq S\cup wz^{-1}S$. Thus,
\begin{eqnarray*}
  P_w\cdot A & = & A\cup wA\\
          & = & w^{-1}S\cup  wz^{-1}S\cup S\cup  w^2z^{-1}S\\
          & = & \{1_H,w^{-1},wz^{-1},w^2z^{-1}\}\cdot S\\
          &=& P_{w^{-1}}\cdot P_{w^2z^{-1}}\cdot S.
\end{eqnarray*}

Applying $F$ to the above equation and using (\ref{3.1}) gives that
\begin{equation}\label{3.3}P_w\cdot B=P_{w^{-1}}\cdot P_{w^2z^{-1}}\cdot T,\end{equation}
where $B=F(A)$. That is,
\begin{equation}\label{3.4}B\cup wB=w^{-1}T\cup T\cup wz^{-1}T\cup w^2z^{-1}T.\end{equation}

By (\ref{3.2}) and (\ref{3.4}), we obtain that
\begin{eqnarray*}
  w^{-2}T\cap B & \subseteq & w^{-2}T\cap(w^{-1}T\cup T\cup wz^{-1}T\cup w^2z^{-1}T)\\
          & = & \cup_{i=-1}^0(w^{-2}T\cap w^i T)\cup(\cup_{j=1}^2(w^{-2}T\cap w^jz^{-1} T))\\
          & = & \emptyset,
\end{eqnarray*}
and thus $w^{-1}T\cap wB=\emptyset$.  So $w^{-1}T\subseteq B$. By the same argument, we have $w^3z^{-1}T\cap wB=\emptyset $ and thus $w^2z^{-1}T\cap B=\emptyset$. It follows that $w^2z^{-1}T\subseteq wB$ and $wz^{-1}T\subseteq B$.

  We have shown $w^{-1}T,wz^{-1}T\subseteq B$ and $w^2z^{-1}T\cap B=\emptyset$.  Combining these with  $(\ref{3.4})$ gives
$$B=w^{-1}T\cup Z\cup wz^{-1}T,$$
where $Z=B\cap T$. Substituting  $B$ into  equation  (\ref{3.4}) yields that
$$wZ\subseteq w^{-1}T\cup T\cup wz^{-1}T\cup w^2z^{-1}T.$$
Since $Z\subseteq T$, we have $wZ\cap w^iT=wZ\cap w^jz^{-1}T=\emptyset$ for any integers $i,j\not=1$ by (\ref{3.2}). We conclude that   $wZ\subseteq wz^{-1}T$ and $Z\subseteq z^{-1}T$. Combing with that  $1_H\in Z$, we obtain $1_H\in z^{-1}T$ and thus $z\in T.$

Case 2: \(z\) is not invertible. By Lemma \ref{lem:3.2} and the same  arguments at the begin of Case 1, we may choose $N>0$ such that
\begin{equation}\label{3.5}w^iR_1\cap w^jR_2=\emptyset,\qquad \text{for}~R_1,R_2\in\{S,T,zS,zT\}~\text{and}~i,j\in\mathbb{Z}~\text{with}~i\neq j,\end{equation}
where $w=g^N$.

Put $$C=\{w^{-2},w^{-1},1_H,w\}.$$ We have the following claim.

{\bf Claim}: For given $R\in\{S,T\}$ and \(X\in\mathcal P_{\mathrm{fin},1}(H)\),
\begin{equation}\label{3.6} P_{w^{-1}}\cdot X=C\cdot R\end{equation}
if and only if
\begin{equation}\label{3.7}X=w^{-1}R\cup V\cup wR,\qquad 1_H\in V\subseteq R.\end{equation}
{\bf Proof of the Claim.}
It is easy to verify that the form given in (\ref{3.7}) satisfies (\ref{3.6}) and belongs to $\mathcal P_{\mathrm{fin},1}(H)$.  Conversely, suppose $A\in\mathcal P_{\mathrm{fin},1}(H)$ and $P_{w^{-1}}\cdot A=C\cdot R$. That is,
\begin{equation}\label{3.8}A\cup w^{-1}A=\cup_{i=-2}^1w^iR.\end{equation}
We need to show  that $w^{-2}R\cap A=\emptyset$, and $w^{-1}R,wR\subseteq A$. If $w^{-2}R\cap A\neq \emptyset$ and $x\in w^{-2}R\cap A$, then $w^{-1}x\in w^{-3}R\cap w^{-1}A$. By (\ref{3.8}), it follows that
$$w^{-1}x\in w^{-3}R\cap(\cup_{i=-2}^{1}w^iR)=\cup_{i=-2}^1(w^{-3}R\cap w^iR),$$ which contradicts to (\ref{3.5}). So $w^{-2}R\cap A=\emptyset$, and thus $w^{-2}R\subseteq w^{-1}A$ by (\ref{3.8}). We obtain $w^{-1}R\subseteq A.$

Similarly, if $wR\cap w^{-1}A\not=\emptyset$, then $wy\in w^2R\cap A$ for any  $y\in wR\cap w^{-1}A$. By (\ref{3.8}), it follows that
$wy\in w^{2}R\cap(\cup_{i=-2}^{1}w^iR)$, which contradicts to (\ref{3.5}). So $wR\cap w^{-1}A=\emptyset$ and we get $wR\subseteq A$ by (\ref{3.8}).

Combing (\ref{3.8}) with that $w^{-2}R\cap A=\emptyset$, and $w^{-1}R,wR\subseteq A$, we have
$$A=\cup_{i=-2}^1(A\cap w^iR)=w^{-1}R\cup V\cup wR,~\qquad 1_H\in V=A\cap R\subseteq R.$$
This completes the proof of the Claim.\\

Let $X=w^{-1}S\cup P_z\cup wS$. Since \(P_z\subseteq S\), by the Claim  we have
\begin{equation}\label{3.9}P_{w^{-1}}\cdot X=C\cdot S.\end{equation}
Note that $C=P_{w^{-1}}\cdot P_{w^{-1}}\cdot P_w$.  Applying $F$ to (\ref{3.9}) and using (\ref{3.1}) gives that
$$P_{w^{-1}}\cdot Y=C\cdot T,$$
where \(Y=F(X)\).
By the Claim  again,
$$Y=w^{-1}T\cup V\cup wT,\qquad 1_H\in V\subseteq T.$$
 By (\ref{3.5}), both $X=w^{-1}S\cup P_z\cup wS$ and $Y=w^{-1}T\cup V\cup wT$ are disjoint unions.
Since \(F\) is cardinality-preserving by Lemma \ref{lemma:2.3}, we assert that $|S|=|T|$ and
$$|X|=2|S|+2=|Y|=2|T|+|V|.$$
It follows that $|V|=2$ and thus $V=P_v$  for some  non-identity element $v\in T$.

We are now in a position to show that $v=z.$ By (\ref{3.5}), we have $R_1\cap zR_2=\emptyset$ for distinct $R_1,R_2\in\{w^{-1}S,P_z,wS\}$. By $(2),(3)$ of Lemma \ref{lemma:3.1},
\begin{eqnarray*}
    \nu_z(X) &=& \nu_z(w^{-1}S\cup P_z\cup wS)  \\
     &=& \nu_{z}(w^{-1}S)+\nu_{z}(P_z)+\nu_{z}(wS) \\
     &=& 2\nu_{z}(S)+\nu_{z}(P_z).
  \end{eqnarray*}
By the same discussion, we also have
$$\nu_z(Y)=2\nu_z(T)+\nu_z(P_v).$$
By Lemma~\ref{lem:3.3}, \(\nu_z(S)=\nu_z(T)\) and
\(\nu_z(X)=\nu_z(Y)\). Hence
$$\nu_z(P_z)=\nu_z(P_v)\in\{1,2\}.$$
Since $z$ is not invertible, $z^2\not\in P_z$ and \(\nu_z(P_z)=1\). Therefore,
$$\nu_z(P_v)=|\{1_H,v\}\cap\{z,zv\}|=1.$$
It is clear that $zv\in\{1_H,v\}$ implies $z$ is invertible. So $zv\not\in P_v$ and $z\in P_v$. We thus get $z=v\in T.$

We have shown that  \(S\subseteq T\). Since \(|S|=|T|<\infty\), we have \(S=T\). Therefore, $F$ is the
identity map.
\end{proof}
\end{theorem}

\begin{theorem}
\label{thm:abelian}
Let \(H\) be an abelian group which is not isomorphic to the Klein four-group. Then
\[
\Phi:\operatorname{Aut}(H)\longrightarrow
\operatorname{Aut}(\mathcal P_{\mathrm{fin},1}(H)),~~~f\longmapsto f^*
\]
is an isomorphism.
\begin{proof}
By Theorem~\ref{lemma:2.4} and Theorem \ref{thm:infinite-unit}, we only need to deal with the case that $H$ is infinite and  torsion.  By Theorem \ref{bothgroup}, the pullback of an automorphism of \(\mathcal P_{\mathrm{fin},1}(H)\) is also an automorphism of $H$ since $H$ is torsion. By Lemma \ref{thm:2.5}, it suffices to prove that if $F$ is an automorphism of $\mathcal P_{\mathrm{fin},1}(H)$  with identity pullback, then $F$ is the identity map.

Let \(S\in\mathcal P_{\mathrm{fin},1}(H)\)  and let $P_h=\{1_H,h\}$ for  $h\in H.$ Since \(H\) is infinite,  abelian, and torsion, there exists a finite subgroup
\(G\) of $H$ such that \(S\subseteq G\) and \(|G|>4\). Since $G=\prod_{g\in G}P_g,$ and $F$ fixes $P_g$, it follows that

 $$F(G)=F(\prod_{g\in G}P_g)=\prod_{g\in G}F(P_g)=\prod_{g\in G}P_g=G.$$ Since $G$ is a subgroup of $H$,   we have $A\cdot G=G$ if and only if $A\subseteq G$  for any $A\in \mathcal{P}_{\mathrm{fin},1}(H)$. Therefore \(F\) maps \(\mathcal P_{\mathrm{fin},1}(G)\) into itself. The same argument applied to \(F^{-1}\) shows that the restriction $F|_{\mathcal P_{\mathrm{fin},1}(G)}$ is also an automorphism of \(\mathcal P_{\mathrm{fin},1}(G)\) with identity pullback. Since \(|G|>4\), by Theorem \ref{lemma:2.4} the restriction $F|_{\mathcal P_{\mathrm{fin},1}(G)}$ itself is the identity map. Hence, $F(S)=S$, which completes the proof.
 \end{proof}
\end{theorem}

Bhowmik and Tringali \cite[Theorem 1.2]{Bho-Tri26} showed that $\mathrm{Aut}(\mathcal{P}_{\mathrm{fin},0}(H))$ is trivial for every numerical monoid $H\not=\mathbb{N}$. Together with our results, this motivates the following conjecture.

\begin{conjecture} Let $H$ be a commutative and cancellative monoid with nontrivial torsion unit group. Suppose that $H$ is not isomorphic to $C_2\times C_2.$ Then all automorphisms of $\mathcal{P}_{\mathrm{fin},1}(H)$ are inner.
\end{conjecture}

\section*{AI Declaration}
During the preparation of this manuscript, the authors used ChatGPT and Kimi to assist with checking the set equation involving four translates in the proof of Theorem~\ref{thm:infinite-unit}, verifying reference formats, and improving English grammar. The main mathematical results and their proofs were developed by the authors independently of these tools. The authors reviewed and edited all AI-assisted material and takes full responsibility for the final manuscript.

\end{document}